%% file: noises_main.tex
\documentclass[a4paper]{article}

\usepackage{xspace,colortbl}
\usepackage{amssymb,amsbsy}
\usepackage{graphicx}
\usepackage[T1]{fontenc}
\usepackage[english]{babel}
\usepackage{makeidx}
\makeindex
\usepackage{mathrsfs}
\usepackage{amsfonts,amsthm}
\usepackage[mathscr]{eucal}
\usepackage{cite}
\usepackage{textcomp}
\usepackage{url}
\theoremstyle{definition}
\newtheorem{ese}{Example}[section]

\newtheorem{rem}{Remark}[section]
\theoremstyle{plain}
\newtheorem{defi}{Definition}[section]
\newtheorem{prop}{Proposition}[section]

\newtheorem{teor}{Theorem}[section]

\begin{document}
\hoffset -2.5cm

\include{noises}

\include{noises_biblio}
\end{document}

%% file: noises.tex
\begin{center}
{\Large \bf A class of self-similar stochastic processes with stationary increments to model anomalous diffusion in physics.}
\end{center}
\vspace{1cm}
\centerline{\bf A Mura$^1$ and F Mainardi$^1$}
\vskip.4cm
\begin{center}
{\it   $1.$ Department of Physics, University of Bologna, and INFN,
Via Irnerio 46, I-40126 Bologna, Italy}
\date{\today}
\end{center}
\vspace{0.5cm}
{\bf Abstract}: 
In this paper we present a general mathematical construction that allows us to define a parametric class of $H$-sssi stochastic processes (self-similar with stationary increments), which have marginal probability density function that evolves in time according to a partial integro-differential equation of fractional type. This construction is based on the theory of finite measures on functional spaces. Since the variance evolves in time as a power function, these $H$-sssi processes naturally provide models for slow and fast anomalous diffusion. Such a class includes, as particular cases, fractional Brownian motion, grey Brownian motion and Brownian motion. 
\vspace{0.5cm}
\section{Introduction}\label{s1}
The  grey noise theory introduced by Schneider (see \cite{Schneider1, Schneider2}) leads naturally to a class of self-similar stochastic processes $\{B_\beta(t),\;0<\beta\le 1\}$. These processes, called {\it grey Brownian motion}, provide stochastic models for the slow-anomalous diffusion\footnote{Anomalous diffusion is characterized by the (asymptotic) time power-law behavior of the variance: $\sigma^2(t)\sim t^\gamma$. Namely, the diffusion is slow if the exponent $\gamma$ is lesser than one, normal if it is equal to one and fast if it is greater than one.} described by the time fractional diffusion equation; i.e.~the marginal density function of the grey Brownian motion is the fundamental solution of the time fractional diffusion equation (see  \cite{Schneider3} and \cite{mainadiitsf,lumapa}). This class will be extended to a class $\{B_{\alpha,\beta}(t)\}$, with $0<\alpha<2,\;\; 0<\beta\le 1$, called {\it ``generalized'' grey Brownian motion}, which includes stochastic models either for slow and fast-anomalous diffusion. First, we present and motivate the mathematical construction. Then, we show that this class is made up of $H$-sssi processes and contain either Gaussian and non-Gaussian processes (like fractional Brownian motion and grey Brownian motion). Finally, we show how the time evolution of the marginal density function is described by partial integro-differential equations of fractional type.\\

We begin introducing some basic concepts and facts. Let $X$ be a vector space over a $\mathbb{K}$-field and let $\left\{||\cdot ||_p,\; p\in I\right\}$ be a countable family of Hilbert-norms  defined on it. The space $X$ along with the Hilbert-norms $\{||\cdot||_p,\;p\in I\}$ is said a {\it topological vector space} if it carries as natural topology the initial topology\footnote{The coarsest topology defined on $X$ which makes these functions continuous.} of the norms and the vector space operations. We indicate with $X_p$ the completion of $X$ with respect to the norm $||\cdot||_p$. Let $\langle \cdot, \cdot \rangle$ denote the natural bilinear pairing between $X$ and its dual space $X'$. We equip $X'$ with the so called {\it weak topology}, which is the coarsest topology such that the functional $\langle \cdot, x\rangle$ is continuous for any $x\in X$.
\begin{defi}[Nuclear space]
A topological vector space $X$, with the topology defined by a family of Hilbert-norms, is said a {\it nuclear space} if for any Hilbert-norm $||\cdot||_p$ there exists a larger norm $||\cdot ||_q$ such that the inclusion map $X_q\hookrightarrow X_p$ is an Hilbert-Schmidt operator\footnote{An Hilbert-Schmidt operator is a bounded operator $A$, defined on an Hilbert space $H$, such that there exists an orthonormal basis $\{e_{i}\}_{i\in I}$ of $H$ with the property $\displaystyle\sum_{i\in I}||Ae_i||^2<\infty$.}.    
\end{defi}
Nuclear spaces have many of the good properties of the finite dimensional Euclidean spaces $\mathbb{R}^d$. For example, a subset of a nuclear space is compact if and only if is bounded and closed. Moreover, spaces whose elements are ``smooth'' is some sense tend to be nuclear spaces.
In the following example we see how nuclear spaces could be constructed naturally starting from an Hilbert space and an operator (see Kuo \cite{kuo}).  
\begin{ese}\label{ese1}
Let $H$ be an Hilbert space and $A$ an operator defined on it. Suppose that there exists an orthonormal bases $\{h_n,\; n=1,2,\dots\}$ satisfying the following properties:
\begin{enumerate}
\item They are eigenvectors of $A$; i.e. for any $n>0$: $Ah_n=\lambda_nh_n, \;\;\lambda_n \in \mathbb{R}$.
\item $\{\lambda_{n}\}_{n>0}$ is a non-decreasing sequence such that: $1\le\lambda_1\le\lambda_2\le \cdots \le \lambda_n$
\item There exists a positive integer $a$ such that:  $\displaystyle\sum_{n=1}^{\infty}\lambda_n^{-a}<\infty$. 
\end{enumerate}
For any non-negative rational number $p\in \mathbb{Q}_+$, we define a sequence of norms $\{||\cdot||_p,\; p\in \mathbb{Q}_+\}$ such that: $||\xi||_p=||A^p\xi||,\;\; \xi \in H$. That is:
\begin{equation}
||\xi||_p=\left(\sum_{n=1}^{\infty}\lambda^{2p}_n(\xi,h_n)^2\right)^{1/2},
\end{equation}
where $(\cdot,\cdot)$ indicates the $H$ inner product. 
\begin{rem}
For any $p\in \mathbb{Q}_+$, the norm $||\cdot ||_p$ is an Hilbert-norm. Indeed, it comes from the scalar product:
\begin{equation}
(\xi,\eta)_p=\sum_{n=1}^{\infty}\lambda_n^{2p}(\xi,h_n)(\eta,h_n).
\end{equation}
\end{rem}
\noindent 
For any $p\in \mathbb{Q}_+$ we define: $X_p=\{\xi \in H;\; ||\xi||_p<\infty\}$. In view of the above remark, $X_p$ is an Hilbert space. Moreover, it is easy to see that for any $p\ge q\ge 0$:
\begin{equation}
X_p\subset X_q.
\end{equation}
We have the following proposition: 
\begin{prop}
For any $p\in \mathbb{Q}_+$, the inclusion map $X_{p+a/2}\hookrightarrow X_p$ is an Hilbert-Schmidt operator.
\end{prop}
\noindent
{\bf Proof}: we set $h^p_n=\displaystyle\frac{1}{\lambda_n^p}h_n$. The collection $\{h_n^p,\; n=1,2,\dots \}$ is an orthonormal bases of $X_p$. In fact, for any positive integers $n$ and $m$:
$$
(h_n^p,h_m^p)_p=\sum_{k=1}^{\infty}\lambda_k^p(h_n^p,h_k)(h_m^p,h_k)=\sum_{k=1}^{\infty}\frac{\lambda_k^{2p}}{\lambda_n^p\lambda_m^p}\delta_{nk}\delta_{mk}=\delta_{nm}.
$$
For each $\xi\in X_{p+a/2}$, we indicate with $i(\xi)=\xi\in X_p$ the inclusion map. Therefore, for any $n>0$:
$$
i(h^{p+a/2}_n)=h^{p+a/2}_n=\frac{1}{\lambda_{n}^{p+a/2}}\lambda_{n}^{p}h^{p}_n=\lambda^{-a/2}h_n^p,
$$
and thus by  hypothesis 
$$
\sum_{n=1}^{\infty}||i(h_n^{p+a/2})||_p^2=\sum_{n=1}^{\infty}\lambda_n^{-a}<\infty,\;\; \Box
$$
Consider the vector space $X=\displaystyle\bigcap_{p\in \mathbb{Q}_+} X_p$. In view of the above proposition $X$ along with the family of Hilbert-norms $\{||\cdot||_p\,,\; p\in \mathbb{Q}_+\}$ is a nuclear space.
\end{ese}
Let $X$ be a vector space. A continuous map $\Phi:X\rightarrow \mathbb{C}$ is called a characteristic functional on $X$ if it's normalized:
$$
\Phi(0)=1,
$$
and positive defined:
$$
\sum_{i,j=1}^m\overline c_i\Phi(\xi_i-\xi_j)c_j\ge 0,\;\;\;m\in \mathbb{Z},\;\; \{c_i\}_{i=1,\dots,m}\in \mathbb{C},\;\;\; \{\xi_i\}_{i=1,\dots,m}\in X.
$$
Let $X=\mathbb{R}^n$. The Bochner theorem \cite{Reed} states that for any characteristic functional $\Phi$ defined on $\mathbb{R}^n$, there exists a unique probability measure $\mu$ defined on $\mathbb{R}^n$, such that 
$$
\displaystyle\int_{\displaystyle\mathbb{R}^n}e^{i(x,\,\xi)}d\mu(x)=\Phi(\xi),\;\, \xi\in \mathbb{R}^n.
$$ 
Let now $X$ be a topological vector space. In the characterization of typical configurations of measures on infinite dimensional spaces the so called Minlos theorem plays a very important role. This theorem is an {\it infinite dimensional} generalization of the Bochner theorem:
\begin{teor}[Minlos theorem]
Let $X$ be a nuclear space. For any characteristic functional $\Phi$ defined on $X$ there exists a unique probability measure $\mu$ defined on the measurable space $(X',\mathcal{B})$, where $\mathcal{B}$ is regarded as the Borel $\sigma$-algebra generated by the weak topology on $X'$, such that: 
\begin{equation}
\int_{X'}e^{i\langle \omega,\xi\rangle}d\mu(\omega)=\Phi(\xi),\;\;\; \xi\in X.
\end{equation}  
\end{teor}
Characteristic functional on Hilbert spaces can be defined starting from completely monotonic functions\footnote{A function $F(t)$ is completely monotone if it is non-negative and possesses derivatives of any order such that:
$
(-1)^k\frac{d^k}{dt^k}F(t)\ge 0,\;\;\; t> 0,\;\; k\in \mathbb{Z}_+=\{0,1,2,\dots\}.
$
}. 
In fact we have the following proposition: 
\begin{prop}\label{p1}
Let $F$ be a completely monotonic function defined on the positive real line. Therefore, there exists a unique characteristic functional $\Phi$, defined on a real separable Hilbert space $H$, such that:
$$
\Phi(\xi)=F(||\xi||^2),\;\; \xi \in H.
$$
\end{prop}
\noindent
This is obvious because completely monotonic functions are associated to non-negative measure defined on the positive real line (see Feller \cite{Feller2}). The converse is also true (see Schneider \cite{Schneider1, Schneider2}).   
\section{White noise}
Consider the Schwartz space $\mathcal{S}(\mathbb{R})$. Equip $\mathcal{S}(\mathbb{R})$ with the usual scalar product:
\begin{equation}
(\xi,\eta)=\int_{\mathbb{R}}dt\xi(t)\eta(t),\;\;\;\xi,\eta \in \mathcal{S}(\mathbb{R}).
\label{sca}
\end{equation}
We indicate the completion of $\mathcal{S}(\mathbb{R})$ with respect to eq.~(\ref{sca}) with $\mathcal{S}_0(\mathbb{R})=\mathscr{L}^2(\mathbb{R})$. We consider the orthonormal system $\{h_n\}_{n\ge 0}$ of the Hermite functions:
\begin{equation}
h_{n}(x)=\frac{1}{\sqrt{(2^nn!\sqrt \pi)}}H_n(x)e^{-x^2/2},
\label{herortho}
\end{equation}
where $H_n(x)=(-1)^ne^{x^2}(d/dx)^ne^{-x^2}$ are the Hermite polynomials of degree $n$. Let $A$ be the ``harmonic oscillator'' operator:
\begin{equation}
A=-\frac{d^2}{dx^2}+x^2+1;
\label{harm}
\end{equation}
$A$ is densely defined on $\mathcal{S}_0(\mathbb{R})$ and the Hermite functions are eigenfunctions of $A$:
$$
Ah_n=\lambda_nh_n=(2n+2)h_n,\;\;\; n=0,1,\dots\;.
$$
We observe that $1\le \lambda_0\le \lambda_1\le\dots \le \lambda_n$ and $\sum_{n }\lambda_n^{-2}<\infty $. We are in the condition of Example \ref{ese1}. Therefore, for any non-negative integer $p$, we can define:
$$
||\xi||_p=||A^p\xi||=\left(\sum_{n=0}^{\infty}(2n+2)^{2p}(\xi,h_n)^2\right)^{1/2},
$$ 
where $||\cdot||$ indicates the $\mathscr{L}^2$ norm. The Schwartz space $\mathcal{S}(\mathbb{R})$ could be then ``reconstructed'' as the {\it projective limit} of the Hilbert spaces $\mathcal{S}_p(\mathbb{R})=\{\xi \in \mathscr{L}^2(\mathbb{R});\;\; ||\xi||_p<\infty\}$. That is:
\begin{equation}
\mathcal{S}(\mathbb{R})=\bigcap_{p\ge 0} \mathcal{S}_p(\mathbb{R}).
\end{equation}
Therefore, the topological Schwartz space, with the topology defined by the $||\cdot||_p$ norms, is a nuclear space. Since $\mathcal{S}(\mathbb{R})$ is a nuclear space, we can apply the Minlos theorem in order to define probability measures on its dual space $\mathcal{S}'(\mathbb{R})$. Consider the positive function $F(t)=e^{-t}$, $t\ge 0$. It is obvious that $F$ is a completely monotone function. Therefore, the functional $\Phi(\xi)=F(||\xi||^2)$, $\xi \in \mathscr{L}^2(\mathbb{R})$, defines a characteristic functional on $\mathcal{S}(\mathbb{R})$. By Minlos theorem, there exists a unique probability measure $\mu$, defined on $(\mathcal{S}'(\mathbb{R}), \mathcal{B})$, such that:
\begin{equation}
\int_{\mathcal{S}'(\mathbb{R})}e^{i\langle \omega, \xi\rangle}d\mu(\omega)=e^{-||\xi||^2},\;\; \xi \in \mathcal{S}(\mathbb{R}).
\end{equation}  
The probability space $(\mathcal{S}'(\mathbb{R}),\mathcal{B},\mu)$ is called {\it white noise space} and the measure $\mu$ is called {\it white noise measure}, or standard Gaussian measure, on $\mathcal{S}'(\mathbb{R})$. \\ 

Consider the generalized stochastic process $X$, defined on the white noise space, such that for each test function $\varphi \in \mathcal{S}(\mathbb{R})$:
\begin{equation}
X(\varphi)(\cdot )=\langle \cdot, \varphi\rangle.
\end{equation}
Clearly, for any $\varphi\in \mathcal{S}(\mathbb{R})$, $X(\varphi)$ is a Gaussian random variable with zero mean and variance $E(X(\varphi)^2)=2||\varphi||^2$. Moreover, for any $\varphi, \phi \in \mathcal{S}(\mathbb{R})$:
\begin{equation}
E(X(\varphi)X(\phi))=2(\varphi,\phi),
\end{equation}
where $E(w)$ indicates the expectation value of the random variable $w$. We refer to the generalized process $X$ as the {\it canonical noise} of $(\mathcal{S}'(\mathbb{R}),\mathcal{B},\mu)$. 
\begin{rem}
In view of the above properties the process $X$ is a white noise \cite{kuo}, and this also motivate the name ``white noise space'' for the probability space $(\mathcal{S}'(\mathbb{R}),\mathcal{B},\mu)$.
\end{rem}
\noindent
We have the following:
\begin{prop} For any $h\in \mathscr{L}^2(\mathbb{R})$, $X(h)$ is defined almost everywhere on $\mathcal{S}'(\mathbb{R})$. Moreover, it is Gaussian with zero mean and variance $2||h||^2$.
\end{prop}
\noindent
{\bf Proof}: we indicate with $(\mathscr{L}^2)=\mathscr{L}^2(\mathcal{S}'(\mathbb{R}),\mu)$. Clearly, for any $\xi \in \mathcal{S}(\mathbb{R})$, we have that $X(\xi)\in (\mathscr{L}^2)$ and:
\begin{equation}
||X(\xi)||^2_{(\mathscr{L}^2)}=E(X(\xi)^2)=2||\xi||^2_{\mathscr{L}^2}.
\label{cau}
\end{equation}
For each $h\in \mathscr{L}^2(\mathbb{R})$, there exists a sequence $\{\xi_n\}_{n\in \mathbb{N}}$ of $\mathcal{S}(\mathbb{R})$-elements which converges to $h$ in the $\mathscr{L}^2(\mathbb{R})$-norm. Therefore, from eq.~(\ref{cau}), the sequence $\{X(\xi_n)\}_{n\in \mathbb{N}}$ is Cauchy in $(\mathscr{L}^2)$ and converges to a limit function $X(h)$, defined on $\mathcal{S}'(\mathbb{R})$. $\Box$\\

The latter proposition states that for every sequence $\{f_t\}_{t\in \mathbb{R}}$ of $\mathscr{L}^2(\mathbb{R})$-functions, depending continuously on a real parameter $t\in \mathbb{R}$, there exists a Gaussian stochastic process
\begin{equation}
\{Y(t)\}_{t\in \mathbb{R}}=\{X(f_t)\}_{t\in \mathbb{R}},
\end{equation}     
defined on the probability space $(\mathcal{S}'(\mathbb{R}),\mathcal{B},\mu)$, which has zero mean, variance $E(Y_t)^2=2||f_t||^2$ and covariance $E(Y(t_1)Y(t_2))=2(f_{t_1},f_{t_2})$.
\begin{rem}
Observe that if $W(x)$, $x\in \mathbb{R}$, is a Wiener process defined on the probability space $(\Omega,\mathcal{F}, P)$, then the functional
\begin{equation}
X(\varphi)=\int \varphi(x)dW(x),\;\; \varphi \in \mathscr{L}^2(\mathbb{R}),
\end{equation}
is a white noise on the space $(\Omega, \mathcal{F}, P)$. Therefore, if we indicate with $1_{[0,t)}(x)$, $t\ge 0$, the indicator function of the interval $[0,t)$, the process 
\begin{equation}
X(1_{[0,t)})=\int_{0}^{t}dW(x)=W(t),\;\; t\ge 0,
\label{noisb}
\end{equation}
is a one-sided Brownian motion.
\end{rem}
\begin{ese}[Brownian motion]
Let $X$ be a white noise defined canonically on the white noise space $(\mathcal{S}'(\mathbb{R}),\mathcal{B},\mu)$. Looking at eq.~(\ref{noisb}), its natural to think that the stochastic process 
\begin{equation}
\{B(t)\}_{t\ge 0}=\{X(1_{[\,0,\,t)})\}_{t\ge 0}, 
\end{equation}
defines a ``standard'' Brownian motion\footnote{With the word ``standard'' Brownian motion we mean that $E(B(1)^2)=2$.}. Indeed, the process $\{X(1_{[\,0,\,t)})\}_{t\ge 0}$, is Gaussian with covariance:
$$
E\left[X(1_{[\,0,\,t)})X(1_{[\,0,\,s)})\right]=2(1_{[\,0,\,t)},1_{[\,0,\,s)})=2\min(t,s),\;\; t,s\ge 0.
$$ 
\end{ese}
\begin{ese}[Fractional Brownian motion]
The stochastic process:
\begin{equation}
\{B_{\alpha/2}(t)\}_{t\ge 0}=\{X(f_{\alpha,\,t})\}_{t\ge 0},\;\;  0<\alpha<2,
\end{equation}
where 
\begin{equation}
f_{\alpha,t}(x)=\frac{1}{C_1(\alpha)}\left((t-x)_+^{\frac{\alpha-1}{2}}-(-x)_+^{\frac{\alpha-1}{2}}\right),\;\; x_+=\max(x,0),
\end{equation}
and  
\begin{equation}
C_1(\alpha)=\frac{\Gamma(\frac{\alpha+1}{2})}{\left(\Gamma(\alpha+1)\sin\frac{\pi\alpha}{2}\right)^{1/2}},
\end{equation}
is a ``standard'' fractional Brownian motion of order $H=\alpha/2$ (see Taqqu \cite{taqqu1}). 
\end{ese} 
\section{Grey noises}
We have seen that white noise is a generalized stochastic process $X$ defined canonically on the white noise space $(\mathcal{S}'(\mathbb{R}),\mathcal{B},\mu)$, with space of test functions $\mathscr{L}^2(\mathbb{R})$. We have remarked that the white noise could also be defined starting from stochastic integrals with respect to the Brownian motion. In this case the space of test function turns out to be the space of integrands of the stochastic integral. Then, the Brownian motion $B(t)$ could be obtained from the white noise by setting $B(t)=X(1_{[0,t)})$. We want to generalize the previous construction in order to define a general class of $H$-sssi processes which includes, Brownian  motion, fractional Brownian motion and more general processes. \\

Consider a one-sided fractional Brownian motion $\{B_{\alpha/2}(t)\}_{t\ge 0}$ with self-similarity parameter $H=\alpha/2$ and $0<\alpha<2$, defined on a certain probability space $(\Omega, \mathcal{F},P)$. The fractional Brownian motion has a spectral representation \cite{taqqu1}:
\begin{equation}
B_{\alpha/2}(t)=\sqrt{C(\alpha)}\int_{\mathbb{R}}\frac{1}{\sqrt{2\pi}}\frac{e^{itx}-1}{ix}|x|^{\frac{1-\alpha}{2}}d\widetilde B(x),\;\; t\ge 0,
\end{equation}  
where $d\widetilde B(x)$ is a complex Gaussian measure such that $d\widetilde B(x)=dB_1(x)+idB_2(x)$ with $dB_1(x)=dB_1(-x)$, $dB_2(x)=-dB_2(-x)$ and where $B_1$ and $B_2$ are independent Brownian motion. Moreover, 
\begin{equation}
C(\alpha)=\Gamma(\alpha+1)\sin\frac{\pi\alpha}{2}.
\label{Calph}
\end{equation}
We observe that 
\begin{equation}
\frac{1}{\sqrt{2\pi}}\frac{e^{itx}-1}{ix}=\widetilde 1_{[0,t)}(x),
\label{fun}
\end{equation}
where we have indicated with $\widetilde f(x)$ the Fourier transform of the function $f$ evaluated on $x\in \mathbb{R}$:
\begin{equation}
\widetilde f(x)=\mathscr{F}\left(f\right)(x)=\frac{1}{\sqrt{2\pi}}\int_{\mathbb{R}}e^{ixy}f(y)dy.
\end{equation}
In view of eq.~(\ref{fun}) we have:
\begin{equation}
B_{\alpha/2}(t)=\sqrt{C(\alpha)}\int_{\mathbb{R}}\widetilde 1_{[0,t)}(x)|x|^{\frac{1-\alpha}{2}}d\widetilde B(x).
\end{equation}
Therefore, if one defines a generalized stochastic process $X$ such that for a suitable choice of a test function $\varphi$
\begin{equation}
X_\alpha(\varphi)=\sqrt{C(\alpha)}\int_\mathbb{R}\widetilde \varphi(x) |x|^{\frac{1-\alpha}{2}}d\widetilde B(x),
\end{equation}
one can write:
\begin{equation}
B_{\alpha/2}(t)=X_\alpha(1_{[0,t)}),\;\; t\ge 0.
\end{equation}
\begin{rem} 
The space of test function can be the space
\begin{equation}
\widetilde \Lambda_\alpha=\{f\in \mathscr{L}^2(\mathbb{R});\;\; ||f||_\alpha^2=C(\alpha)\int_{\mathbb{R}}|\widetilde f(x)|^2|x|^{1-\alpha}dx<\infty\},
\label{normal}
\end{equation} 
which coincides with a space of deterministic integrands for  fractional Brownian motion (see Pipiras and Taqqu \cite{PipirasTaqqu1, PipirasTaqqu2}).
\end{rem}
\noindent
Consider now the Schwartz space $\mathcal{S}(\mathbb{R})$  equipped with the scalar product:
\begin{equation}
(\xi,\eta)_\alpha=C(\alpha)\int_{\mathbb{R}}\overline{\widetilde\xi(x)}\widetilde\eta(x)|x|^{1-\alpha}dx,\;\;\;\xi,\eta \in \mathcal{S}(\mathbb{R}),\;\; 0<\alpha<2,
\label{scaalpha}
\end{equation}
where $C(\alpha)$ is given by eq.~(\ref{Calph}). This scalar product generate the $\alpha$-norm in eq.~(\ref{normal}). We indicate with $\mathcal{S}_0^{(\alpha)}(\mathbb{R})$ the completion of $\mathcal{S}(\mathbb{R})$ with respect to eq.~(\ref{scaalpha}).
\begin{rem} \label{r1}
If we set $\alpha=1$ in eq.~(\ref{scaalpha}), we have $C(1)=1$ and:
\begin{equation}
(\xi,\eta)_1=\int_{\mathbb{R}}\overline{\widetilde\xi(x)}\widetilde\eta(x)dx=\int_{\mathbb{R}}\xi(y)\eta(y)dy,
\end{equation}
so that, we recover the $\mathscr{L}^2(\mathbb{R})$-inner product. Moreover, $\mathcal{S}_0^{(1)}(\mathbb{R})=\mathcal{S}_0(\mathbb{R})=\mathscr{L}^2(\mathbb{R})$.
\end{rem}
\noindent
Starting from the Hilbert space $(\mathcal{S}_0^{(\alpha)}(\mathbb{R}),||\cdot||_\alpha)$, it is possible to reproduce the construction of Example \ref{ese1}. Then, the space $\mathcal{S}(\mathbb{R})$ turns out to be a  nuclear space with respect to the topology generated by the $\alpha$-norm $||\cdot||_\alpha$ and an operator $A^{(\alpha)}$. Here we just say that the main ingredient are the {\it Generalized Laguerre polynomials}:
\begin{equation}
L^{\gamma}_n(x)=\frac{x^{-\gamma}e^{x}}{\Gamma(n+1)}\frac{d^n}{dx^n}\left(e^{-x}x^{n+\gamma}\right),\;\; \gamma>-1,\;\;x \ge 0,
\end{equation}
where $n$ is a non-negative integer. They are orthogonal with respect to the weighting function $x^\gamma e^{-x}$,
\begin{equation}
\int_{0}^{\infty} x^{\gamma}e^{-x}L^{\gamma}_n(x)L^{\gamma}_{m}(x)dx=\frac{\Gamma(n+\gamma+1)}{\Gamma(n+1)}\delta_{nm},
\label{orthola}
\end{equation}
and satisfy the Laguerre equation:
\begin{equation}
\left(x\frac{d^2}{dx^2}+(\gamma+1-x)\frac{d}{dx}\right)L^{\gamma}_n(x)=-nL_n^{\gamma}(x).
\label{laguerre}
\end{equation}
Using eq.~(\ref{orthola}), it is easy to show that the sequence of functions $\{h^{\alpha}_{n}\}_{n\in \mathbb{Z}_+}$ defined by:
\begin{equation}
\left\{
\begin{array}{ll}
\widetilde h^{\alpha}_{2n}(x)=a_{n,\alpha}e^{-x^2/2}L^{-\alpha/2}_n(x^2), & n\in \mathbb{Z}_+;\\[0.3cm]
\widetilde h^{\alpha}_{2n+1}(x)=b_{n,\alpha}e^{-x^2/2}xL^{1-\alpha/2}_n(x^2),& n\in \mathbb{Z}_+,
\end{array}
\right.
\label{tred}
\end{equation}
is an orthonormal bases of $\mathcal{S}_0^{(\alpha)}(\mathbb{R})$ with the choice:
\begin{equation}
a_{\alpha,n}=\left(\frac{\Gamma(n+1)}{C(\alpha)\Gamma(n+1-\alpha/2)}\right)^{1/2},\;\;\;b_{\alpha,n}=\left(\frac{\Gamma(n+1)}{C(\alpha)\Gamma(n+2-\alpha/2)}\right)^{1/2}.
\end{equation}
Then, using eq.~(\ref{laguerre}), one can show that the orthonormal bases $\{h^{\alpha}_{n}\}_{n\in \mathbb{Z}_+}$ is a set of eigenfunction of an operator $A^{(\alpha)}$, defined on $\mathcal{S}_0^{(\alpha)}(\mathbb{R})$, with eigenvalues $\lambda_n^{(\alpha)}=2n+2-\alpha+1$.  
\begin{rem}\label{r2}
We recall the well known relationships between Laguerre and Hermite polynomials:
\begin{equation}
\left\{
\begin{array}{ll}
H_{2n}(x)=(-1)^n2^{2n}n!L_{n}^{-1/2}(x^2)\\[0.3cm]
H_{2n+1}(x)=(-1)^n2^{2n+1}n!xL_{n}^{1/2}(x^2).
\end{array}
\right.
\label{lahe}
\end{equation}
In view of the above relations, when $\alpha=1$ the orthonormal bases $\{h^{\alpha}_{n}\}_{n\in \mathbb{Z}_+}$ reduces to the Hermite bases of $\mathscr{L}^2(\mathbb{R})$ eq.~(\ref{herortho}), which is preserved under Fourier transformation.     
\end{rem}
\noindent
By Proposition \ref{p1}, starting from a completely monotonic function $F$, we can define characteristic functionals on $\mathcal{S}(\mathbb{R})$ by setting $\Phi(\xi)=F(||\xi||^2_\alpha)$. Then, we could use Minlos theorem in order to define probability measures on $\mathcal{S}'(\mathbb{R})$. We consider the real valued Mittag-Leffler function of order $\beta>0$:
\begin{equation}
E_\beta(x)=\sum_{n=0}^{\infty}\frac{x^n}{\Gamma(\beta n+1)},\;\; x\in \mathbb{R}.
\label{equ1}
\end{equation}
It is known that the function $F_\beta(t)=E_\beta(-t)$, $t\ge 0$, is a completely monotonic function if $0<\beta\le 1$ \cite{miller}. For example if $\beta=1$ we recover $F_{1}(t)=e^{-t}$. Therefore, the functional $\Phi_{\alpha,\beta}(\xi)=F_\beta(||\xi||^2_\alpha)$, $\xi \in \mathcal{S}_{0}^{(\alpha)}(\mathbb{R})$, defines a characteristic functional on $\mathcal{S}(\mathbb{R})$. By Minlos theorem, there exists a unique probability measure $\mu_{\alpha,\beta}$, defined on $(\mathcal{S}'(\mathbb{R}),\mathcal{B})$, such that:
\begin{equation}
\int_{\mathcal{S}'(\mathbb{R})}e^{i\langle \omega,\xi\rangle}d\mu_{\alpha,\beta}(\omega)=F_{\beta}(||\xi||_\alpha^2),\;\; \xi\in \mathcal{S}(\mathbb{R}).
\label{gmes}
\end{equation} 
When $\alpha=\beta$ and $0<\beta\le 1$, the probability space $(\mathcal{S}'(\mathbb{R}), \mathcal{B},\mu_{\beta,\beta})$ is called {\it grey noise space} and the measure $\mu_{\beta,\beta}$ is called {\it grey noise measure} (see Schneider \cite{Schneider1,Schneider2}). In this paper, we focus on the more general case $0<\alpha<2$ and we call the space  $(\mathcal{S}'(\mathbb{R}), \mathcal{B},\mu_{\alpha,\beta})$  ``generalized'' grey noise space and $\mu_{\alpha,\beta}$ ``generalized'' grey noise measure.
\begin{defi}
The generalized stochastic process $X_{\alpha,\beta}$, defined canonically on the ``generalized'' grey noise space $(\mathcal{S}'(\mathbb{R}), \mathcal{B},\mu_{\alpha,\beta})$, is called ``generalized'' grey noise. Therefore, for each test function $\varphi\in \mathcal{S}(\mathbb{R})$:
\begin{equation}
X_{\alpha,\beta}(\varphi)(\cdot )=\langle \cdot,\varphi\rangle. 
\end{equation}
\end{defi}
\begin{rem}
By the definition of ``generalized'' grey noise measure eq.~(\ref{gmes}), for any $\varphi\in \mathcal{S}(\mathbb{R})$, we have:
\begin{equation}
E(e^{iyX_{\alpha,\beta}(\varphi)})=E_{\beta}(-y^2||\varphi||_{\alpha}^2),\;\; y\in \mathbb{R}.
\label{gnf}
\end{equation}
\end{rem}
\noindent
Using eq.~(\ref{gnf}) and eq.~(\ref{equ1}) it easy to show that the ``generalized'' grey noise has moments of any order:
\begin{equation}
\left\{
\begin{array}{ll}
E(X_{\alpha,\beta}(\xi)^{2n+1})=0, \\[0.3cm]
E(X_{\alpha,\beta}(\xi)^{2n})=\displaystyle\frac{2n!}{\Gamma(\beta n+1)}||\xi||_\alpha^{2n},
\end{array}
\right.
\label{mom}
\end{equation}
for any integer $n\ge 0$ and $\xi\in \mathcal{S}(\mathbb{R})$. It is possible to extend the space of test functions to the whole $\mathcal{S}_{0}^{(\alpha)}(\mathbb{R})$. In fact, for any $\xi \in \mathcal{S}(\mathbb{R})$ we have $X_{\alpha,\beta}(\xi)\in (\mathscr{L}^2)=\mathscr{L}^2(\mathcal{S}'(\mathbb{R}),\mu_{\alpha,\beta})$. Thus, for any $h\in \mathcal{S}_{0}^{(\alpha)}(\mathbb{R})$, the function $X_{\alpha,\beta}(h)$ is defined as a limit of a sequence $X_{\alpha,\beta}(\xi_n)$, where $\{\xi_n\}$ belong to $\mathcal{S}(\mathbb{R})$. Therefore we have the following:
\begin{prop} 
For any $h\in \mathcal{S}_{0}^{(\alpha)}(\mathbb{R})$, $X_{\alpha,\beta}(h)$ is defined almost everywhere on $\mathcal{S}'(\mathbb{R})$ and belongs to $(\mathscr{L}^2)$. 
\end{prop}   
Summarizing: the ``generalized'' grey noise is defined canonically on the grey noise space $(\mathcal{S}'(\mathbb{R}),\mathcal{B},\mu_{\alpha,\beta})$ with the following properties:
\begin{enumerate}
\item for any $h\in \mathcal{S}_0^{(\alpha)}(\mathbb{R})$, $X_{\alpha,\beta}(h)$ is well defined and belong to $(\mathscr{L}^2)$.
\item $E(e^{iyX_{\alpha,\beta}(h)})=E_{\beta}(-y^2||h||_{\alpha}^2)$ for any $y\in \mathbb{R}$.
\item $E(X_{\alpha,\beta}(h))=0$ and $E(X_{\alpha,\beta}(h)^2)=\displaystyle\frac{2}{\Gamma(\beta+1)}||h||_{\alpha}^2$.
\item For any $h$ and $g$ which belong to $\mathcal{S}_0^{(\alpha)}(\mathbb{R})$, one has:
\begin{equation}
E\left(X_{\alpha,\beta}(h)X_{\alpha,\beta}(g)\right)=\frac{1}{\Gamma(\beta+1)}\left[(h,g)_\alpha+\overline{(h,g)}_\alpha\right].
\label{gnaut}
\end{equation}
\end{enumerate}
If we put $\beta=1$, the measure $\mu_{\alpha,1}:=\mu_\alpha$ is a Gaussian measure and $X_{\alpha,1}:=X_{\alpha}$ is a Gaussian noise. In fact, for any $h\in \mathcal{S}_0^{(\alpha)}(\mathbb{R})$, the random variable $X_\alpha(h)$  is Gaussian with zero mean and variance $E(X_\alpha(h)^2)=2||h||_\alpha^2$ (see eq.~\ref{gnf}). When $\alpha=1$, $X_\alpha$ reduces to a ``standard'' white noise (see Remark \ref{r1} and Remark \ref{r2}). Moreover, for any sequence $\{f_t\}_{t\in \mathbb{R}}$ of $\mathcal{S}_0^{(\alpha)}(\mathbb{R})$-functions, depending continuously on a real parameter $t\in \mathbb{R}$, the stochastic process $Y(t)=X_\alpha(f_t)$ is Gaussian with auto-covariance given by eq.~(\ref{gnaut})
\begin{equation} 
E(Y(t)Y(s))=E(X_\alpha)(f_t)X_{\alpha}(f_s)=(f_t,f_s)_{\alpha}+\overline{(f_t,f_s)}_\alpha.
\end{equation}  
\begin{ese}[Fractional Brownian motion] For any $t\ge 0$ the function $1_{[0,t)}$ belongs to $\mathcal{S}_0^{(\alpha)}(\mathbb{R})$. In fact, it is easy to show that $||1_{[0,t)}||_\alpha^2<\infty$ when $0<\alpha<2$ and 
\begin{equation}
||1_{[0,t)}||_\alpha^2=\frac{C(\alpha)}{2\pi}\int_{\mathbb{R}}dx\frac{2}{|x|^{1+\alpha}}(1-\cos tx) =t^\alpha.
\label{indt}
\end{equation} 
Therefore, we can define the process:
\begin{equation}
B_{\alpha/2}(t)=X_\alpha(1_{[0,t)}),\;\; t\ge 0.
\end{equation} 
The process $B_{\alpha/2}(t)$ is a ``standard'' fractional Brownian motion with parameter $H=\alpha/2$. Indeed, it is Gaussian with variance $E(B_{\alpha/2}(t)^2)=2||1_{[0,t)}||_\alpha^2=2t^\alpha$ and auto-covariance:
$$
E(B_{\alpha/2}(t)B_{\alpha/2}(s))=(1_{[0,t)},1_{[0,s)})_{\alpha}+\overline{(1_{[0,t)},1_{[0,s)})}_\alpha
$$
$$
=\frac{C(\alpha)}{2\pi}\int_{\mathbb{R}}dx \frac{2}{|x|^{\alpha+1}}\left(1-\cos tx+1-\cos sx-1+\cos(t-s)x\right)
$$
$$
=t^\alpha+s^\alpha-|t-s|^\alpha=\gamma_\alpha(t,s),\;\; t,s\ge 0,
$$
which is the fractional Brownian motion auto-covariance.  
\end{ese}
\noindent
In view of the above example, $X_\alpha$ could be regarded as a {\it fractional Gaussian noise} defined on the space $(\mathcal{S}'(\mathbb{R}),\mathcal{B},\mu_\alpha)$.
\begin{ese}[Deconvolution of Brownian motion]
The stochastic process 
\begin{equation}
\{B(t)\}_{t\ge 0}=\{X_\alpha(g_{\alpha,t})\}_{t\ge 0},
\label{dec}
\end{equation}
where, for each $t\ge 0$, the function $g_{\alpha,t}$ is defined by: 
\begin{equation}
\widetilde g_{\alpha,t}(x)=\frac{1}{\sqrt{C(\alpha)}}\widetilde 1_{[0,t)}(x)(ix)^{\frac{\alpha-1}{2}}, 
\end{equation}
is a ``standard'' Brownian motion. Indeed, it is Gaussian, with zero mean, variance
\begin{equation}
E(B(t)^2)=2\int_{\mathbb{R}}|x|^{1-\alpha}|\widetilde 1_{[\,0,t)}(x)|^2|x|^{\alpha-1}dx=2\int_{\mathbb{R}}|\widetilde 1_{[\,0,t)}(x)|^2dx=2t,
\end{equation}
and autocovariance:
\begin{equation}
E(B(t)B(s))=\int_{\mathbb{R}}\left(\overline{\widetilde 1_{[0,t)}}(x)\widetilde 1_{[0,s)}(x) +\overline{\widetilde 1_{[0,s)}}(x)\widetilde 1_{[0,t)}(x)\right)dx=2\min(t,s).
\end{equation}
\end{ese}
\begin{rem}
The representation of Brownian motion in terms of the fractional Gaussian noise eq.~(\ref{dec}) corresponds to a particular case of the so called {\it deconvolution formula}, which expresses the Brownian motion as a stochastic integral with respect to a fractional Brownian motion of order $H=\alpha/2$ (see  \cite{PipirasTaqqu3}). More generally, we can represent a fractional Brownian motion $B_{\gamma/2}(t)$ of order $H=\gamma/2$, $0<\gamma<2$ in terms of a fractional Gaussian noise of order $\alpha$, which corresponds to a representation of $B_{\gamma/2}$ in terms of a stochastic integral of a fractional Brownian motion $B_{\alpha/2}$ of order $H=\alpha/2$, $0<\alpha<2$ (see example below).   
\end{rem}
\begin{ese}[Deconvolution of fractional Brownian motion]
The stochastic process,
\begin{equation}
\{B_{\gamma/2}(t)\}_{t\ge 0}=\{X_\alpha(g_{\alpha,\gamma,t})\}_{t\ge 0},
\end{equation}
where:
\begin{equation}
\widetilde g_{\alpha,\gamma,t}(x)=\sqrt{\frac{C(\gamma)}{C(\alpha)}}\widetilde 1_{[0,t)}(x)(ix)^{\frac{\alpha-\gamma}{2}},\;\; 0<\gamma<2,
\end{equation}
is a ``standard'' fractional Brownian motion of order $H=\gamma/2$.
\end{ese}
We consider now the general case $0<\alpha<2$, $0<\beta\le 1$. 
\begin{defi}
The stochastic process
\begin{equation}
\{B_{\alpha,\beta}(t)\}_{t\ge 0}=\{X_{\alpha,\beta}(1_{[0,t)})\}_{t\ge 0},
\end{equation}  
is called ``generalized'' (standard) grey Brownian motion.
\end{defi} 
\noindent
The ``generalized'' grey Brownian motion $B_{\alpha,\beta}$ has the following properties which come directly from the grey noise properties and eq. (\ref{indt}):
\begin{enumerate}
\item $B_{\alpha,\beta}(0)=0$ almost surely. Moreover, for each $t\ge 0$, $E(B_{\alpha,\beta}(t))=0$  and
\begin{equation}
E(B_{\alpha,\beta}(t)^2)=\displaystyle\frac{2}{\Gamma(\beta+1)}t^\alpha.
\label{varg}
\end{equation}
\item The auto-covariance function is:
\begin{equation}
E(B_{\alpha,\beta}(t)B_{\alpha,\beta}(s))=\gamma_{\alpha,\beta}(t,s)=\frac{1}{\Gamma(\beta+1)}\left(t^\alpha+s^\alpha-|t-s|^\alpha\right).
\label{depd}
\end{equation}
\item For any $t,s \ge 0$, the characteristic function of the increments is:
\begin{equation}
E\left(e^{iy(B_{\alpha,\beta}(t)-B_{\alpha,\beta}(s))}\right)=E_{\beta}(-y^2|t-s|^\alpha),\;\; y\in \mathbb{R}.
\label{chinc}
\end{equation} 

\end{enumerate}
The third property follows from the linearity of the grey noise definition. In fact, suppose $0\le s<t$, we have $y(B_{\alpha,\beta}(t)-B_{\alpha,\beta}(s))=yX_{\alpha,\beta}(1_{[0,t)}-1_{[0,s)})=X_{\alpha,\beta}(y1_{[s,t)})$, and $||y1_{[s,t)}||_\alpha^2=y^2(t-s)^\alpha$. All these properties are enclosed in the following:
\begin{prop}\label{p3}
For any $0<\alpha<2$ and $0<\beta\le 1$, the process $B_{\alpha,\beta}(t)$, $t\ge 0$, is a self-similar with stationary increments process ($H$-sssi), with $H=\alpha/2$.
\end{prop}
\noindent
\begin{figure}[t!]
\begin{center} 
\includegraphics[keepaspectratio=true,height=9cm]{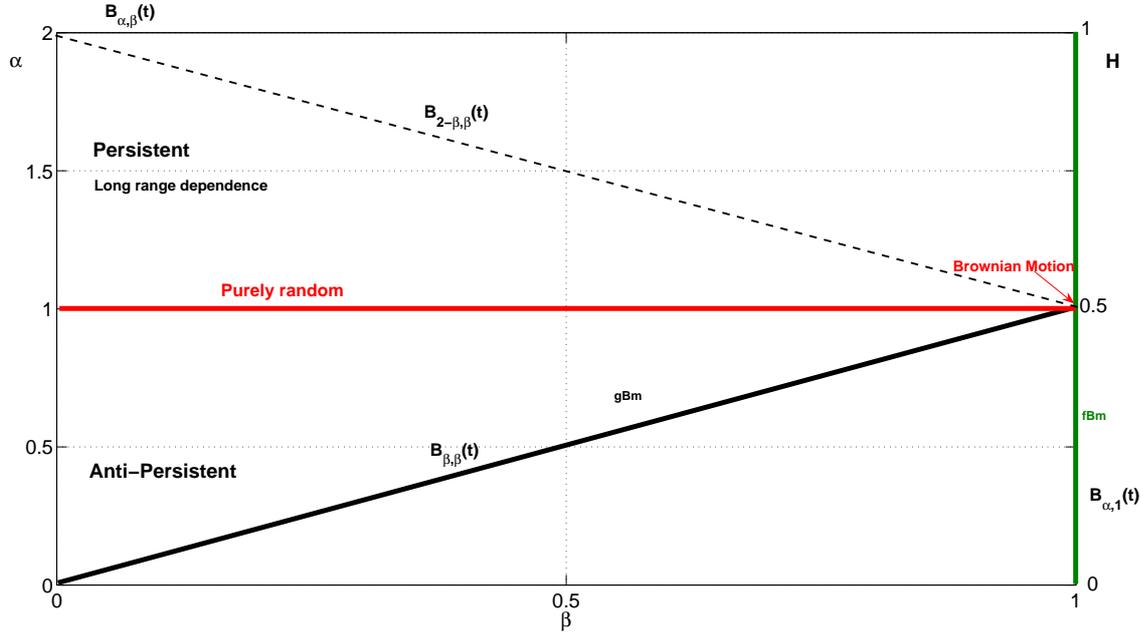}
\end{center}
\caption{Parametric class of generalized grey Brownian motion. The upper diagonal line indicates the ``conjugated'' process of grey Brownian motion.\label{figura 1}} 
\end{figure}
\noindent 
{\bf Proof}: This result is actually a consequence of the linearity of the noise definition.
Given a sequence of real numbers $\{\theta_j\}_{j=1,\dots , n}$, we have to show that for any $0<t_1<t_2<\dots <t_n$ and $a>0$:
$$
E\left(\exp(i\displaystyle\sum_j\theta _jB_{\alpha,\beta}(at_j))\right)=E\left(\exp(i\displaystyle\sum_j\theta _ja^{\frac{\alpha}{2}}B_{\alpha,\beta}(t_j))\right).
$$ 
The linearity of the grey noise definition allows to write the above equality as:
$$
E\left[\exp\left(i X_{\alpha,\beta}\big{(}\sum_j\theta _j 1_{[0,at_j)}\big{)}\right)\right]=E\left[\exp\left(i X_{\alpha,\beta}\big{(}a^{\frac{\alpha}{2}}\sum_j\theta _j 1_{[0,t_j)}\big{)}\right)\right].
$$
Using eq. (\ref{gnf}) we have 
$$
F_{\beta}\left(||\sum_j\theta _j 1_{[0,at_j)} ||_\alpha^2\right)=F_{\beta}\left(||a^{\frac{\alpha}{2}}\sum_j\theta _j 1_{[0,t_j)} ||_\alpha^2\right)
$$
which, because the complete monotonicity, reduces to
$$
||\sum_j\theta _j 1_{[0,at_j)} ||_\alpha^2=a^{\alpha}||\sum_j\theta _j 1_{[0,t_j)} ||_\alpha^2.
$$
In view of the definition eq.~(\ref{scaalpha}) and eq.~(\ref{fun}), the above equality is checked after a simple change of variable in the integration. In the same way we can prove the stationarity of the increments. We have to show that for any $h\in \mathbb{R}$:
$$
E\left[\exp\left(i\displaystyle\sum_j\theta _j(B_{\alpha,\beta}(t_j+h)-B_{\alpha,\beta}(h))\right)\right]=E\left[\exp\left(i\displaystyle\sum_j\theta _j(B_{\alpha,\beta}(t_j)\right)\right].
$$ 
We use the linearity property to write: 
$$
E\left[\exp\left(i X_{\alpha,\beta}\big{(}\sum_j\theta _j 1_{[h,t_j+h)}\big{)}\right)\right]=E\left[\exp\left(i X_{\alpha,\beta}\big{(} \sum_j\theta _j 1_{[0,t_j)}\big{)}\right)\right].
$$
By using the definition and the complete monotonicity, we have:
$$
||\sum_j\theta _j 1_{[h,t_j+h)} ||_\alpha^2=||\sum_j\theta _j 1_{[0,t_j)} ||_\alpha^2
$$
which is true because:
$$
\widetilde 1_{[h,t_j+h)}(x)=\frac{1}{\sqrt{2\pi}}\displaystyle\frac{e^{ixh}}{ix}\left(e^{ixt_j}-1\right).\;\;\; \Box
$$

In view of Proposition \ref{p3}, $\{B_{\alpha,\beta}(t)\}$ forms a class of $H$-sssi stochastic processes indexed by two parameters $0<\alpha<2$ and $0<\beta\le 1$. This class includes fractional Brownian motion ($\beta=1$), grey Brownian motion ($\alpha=\beta$) and Brownian motion ($\alpha=\beta=1$). In Figure \ref{figura 1} we present a diagram which allows us to identify the elements of the class. The long-range dependence\footnote{An $H$-sssi process is said to possess long-range dependence if the discrete process of the increments exhibits long-range dependence. That is, if the increments autocorrelation function tends to zero like a power function and such that it doesn't result integrable \cite{taqqu1}.} domain corresponds to the region $1<\alpha<2$. The horizontal line represents the processes with purely random increments, that is, processes which possess uncorrelated increments. The fractional Brownian motion is identified by the vertical  line ($\beta=1$). The lower diagonal line represents the grey Brownian motion.
\section{Master equation and concluding remarks}
The following proposition characterizes the marginal density function of the process $\{B_{\alpha,\beta}(t), \;t\ge 0\}$: 
\begin{prop}\label{p4}
The marginal probability density function $f_{\alpha,\beta}(x,t)$ of the process $\{B_{\alpha,\beta}(t),\;t\ge 0\}$ is the fundamental solution of the ``stretched'' time-fractional diffusion equation:
\begin{equation}
u(x,t)=u_0(x)+\frac{1}{\Gamma(\beta)}\int_{0}^t\frac{\alpha}{\beta}s^{\alpha/\beta-1}\left(t^{\frac{\alpha}{\beta}}-s^{\frac{\alpha}{\beta}}\right)^{\beta-1}\frac{\partial^2}{\partial x^2}u(x,s)ds, \;\; t\ge 0.
\label{geneq}
\end{equation} 
\end{prop}  
\noindent
{\bf Proof:} eq.~(\ref{chinc}) (with $s=0$) states that $\widetilde f_{\alpha,\beta}(y,t)=E_{\beta}(-y^2t^\alpha)$. Using eq.~(\ref{equ1}), we can show that the Mittag-Leffler function satisfies
$$
E_{\beta}(-y^2(t^{\frac{\alpha}{\beta}})^\beta)=1-\frac{y^2}{\Gamma(\beta)}\int_{0}^{t^{\frac{\alpha}{\beta}}}ds'(t^{\frac{\alpha}{\beta}}-s')^{\beta-1}E_{\beta}(-y^2s'^\beta)
$$
$$
=1-\frac{y^2}{\Gamma(\beta)}\int_{0}^{t}\frac{\alpha}{\beta}s^{\alpha/\beta-1}(t^{\frac{\alpha}{\beta}}-s^{\frac{\alpha}{\beta}})^{\beta-1}E_{\beta}(-y^2s^\alpha)ds,
$$
where we have used the change of variables $s'=s^{\alpha/\beta}$. Thus, $f_{\alpha,\beta}(x,t)$ solves eq.~(\ref{geneq}) with initial condition $u_0(x)=f_{\alpha,\beta}(x,0)=\delta(x)$. $\Box$\\

We refer to eq.~(\ref{geneq}) as the master equation of the marginal density function of the ``generalized'' grey Brownian motion. Therefore, the diagram in Figure \ref{figura 1} can be also read in terms of partial integro-differential equation of fractional type. When $\alpha=\beta$ and $0<\beta\le 1$, we recover the time-fractional diffusion equation of order $\beta$ (lower diagonal line). When $\beta=1$ and $0<\alpha<2$, we have the equation of the fractional Brownian motion marginal density, that is the equation of a stretched Gaussian density (vertical line). Finally, when $\alpha=\beta=1$ we find the standard diffusion equation. \\ 

From Proposition \ref{p4} it follows that the parametric class $\{B_{\alpha,\beta}(t)\}$  provides stochastic models for anomalous diffusions described by eq.~(\ref{geneq}). Looking at eq.~(\ref{varg}) and eq.~(\ref{depd}), which describe the variance and the covariance function respectively, it follows that: 

When $0<\alpha<1$, the diffusion is slow. The increments of the process $B_{\alpha,\beta}(t)$ turn out to be negatively correlated. This implies that the trajectories are very ``zigzaging'' (antipersistent). The increments form a stationary process which does not exhibit long-range dependence. 

When $\alpha=1$, the diffusion is normal. The increments of the process are uncorrelated. The trajectories are said to be ``chaotic''. 

When $1<\alpha<2$, the diffusion is fast. The increments of the process $B_{\alpha,\beta}(t)$ are positively correlated. So that, the trajectories are more regular (persistent). In this case the increments exhibits long-range dependence \cite{taqqu1}.\\

The stochastic processes considered so far, governed by the master equation (\ref{geneq}),  are of course Non-Markovian. We observe that non-Markovian equations like eq.~(\ref{geneq}) are often associated to subordinated stochastic processes $D(t)=B(l(t))$, where the parent Markov process $B(t)$ is a ``standard'' Brownian motion and the random time process $l(t)$ is a self-similar of order $H=\beta$ non-negative non-decreasing  non-Markovian process. For example, in Kolsrud \cite{kolsrud} the random time $l(t)$ is taken to be related to the {\it local time} of a $d=2(1-\beta)$-dimesional fractional Bessel process, while in Meerschaert et al. \cite{meerschaert} (see also Gorenflo et al. \cite{gorenflo} and Stanislavsky \cite{stani}), in the context of Continuous Time Random Walk, it is interpreted as the {\it inverse process} of the totally skewed strictly $\beta$-stable process. Heuristically, our stochastic process $\{B_{\alpha,\beta}(t),\;t\ge 0\}$ cannot be a subordinated process (for example if $\beta=1$ it reduces to a fractional Brownian motion). Therefore, here we provided an example of a class of stochastic models associated to time-fractional diffusion equations like eq.~(\ref{geneq}), which are not subordinated processes. 

It is important to remark that, starting from a master equation which describes the dynamic evolution of a probability density function $f(x,t)$, it is always possible to define an equivalence class of stochastic processes with the same marginal density function $f(x,t)$. All these processes provide suitable stochastic models for the starting equation. In this paper we focused on a subclass $\{B_{\alpha,\beta}(t),\; t\ge 0\}$ associated to the non-Markovian equation eq.~(\ref{geneq}). This subclass is made up of processes with {\it stationary increments}. In this case, the memory effects are enclosed in the typical dependence structure of a $H$-sssi process  eq.~(\ref{depd}); while, for instance in the case of a subordinated process, these are due to the non-Markovian property of the {\it random time process}.

It is also interesting to observe that the ``generalized'' grey Brownian motion turns out to be a direct generalization of a Gaussian process. Indeed, it includes the fractional Brownian motion as particular case when $\beta=1$. Moreover, for any sequence of real numbers $\{\theta_{i}\}_{i=1,\dots, n}$, if one considers the collection $\{B_{\alpha,\beta}(t_1),\dots,B_{\alpha,\beta}(t_n)\}$ with $0<t_1<t_2<\cdots <t_n$, it is easy to show that:
$$
E\left(\exp\Big{(}i\sum_{j=1}^{n}\theta_jB_{\alpha,\beta}(t_j)\Big{)}\right)=E\left(\exp\left(iX_{\alpha,\beta}\Big{(}\sum_{j=1}^{n}\theta_j1_{[0,t_j)}\Big{)}\right)\right)
$$
\begin{equation}
=E_{\beta}\left(-\Big{|}\Big{|}\sum_{j}\theta_j1_{[0,t_j)}\Big{|}\Big{|}_\alpha^2\right)=E_{\beta}\left(-\Gamma(\beta+1)\frac{1}{2}\sum_{i,j}\theta_i\theta_j\gamma_{\alpha,\beta}(t_i,t_j)\right),
\label{furloc}
\end{equation}
where $\gamma_{\alpha,\beta}$ is the the autocovariance matrix eq.~(\ref{depd}). It is clear that, fixed $\beta$, the ``generalized'' grey Brownian motion is defined only by its covariance structure. In other words, $B_{\alpha, \beta}(t)$ provides an example of a stochastic process characterized only by the first and second moments, which is a property of Gaussian processes.